\documentclass[11pt]{article}

\usepackage{amssymb,latexsym}
\usepackage{amsmath,amscd}
\usepackage{theorem}
\usepackage{eucal}

\title{Hermitian star products are completely positive deformations}

\author{\textbf{Henrique Bursztyn}\thanks{E-mail:
    henrique@math.toronto.edu}\addtocounter{footnote}{5}
  \\[0.1cm]
  Department of Mathematics\\
  University of Toronto\\
  100 St. George Street\\
  Toronto, Ontario, M5S 3G3\\[0.1cm]
  and\\[0.1cm]
  The Fields Institute\\
  222 College St.\\
  Toronto, Ontario, MST 3J1\\ 
  Canada\\
[0.5cm]
  \textbf{Stefan Waldmann}\thanks{E-mail:
    Stefan.Waldmann@physik.uni-freiburg.de}
  \\[0.1cm]
  Fakult{\"a}t f{\"u}r Mathematik und Physik\\
  Albert-Ludwigs-Universit{\"a}t Freiburg\\
  Physikalisches Institut\\
  Hermann Herder Stra{\ss}e 3\\
  D 79104 Freiburg\\
  Germany}

\date{October 2004\\[0.5cm] FR-THEP 2004/18}

\newcommand{\textdef}[1] {\textbf{#1}}

\newcommand{\im}         {{\mathrm i}}
\newcommand{\eu}         {{\mathrm e}}

\newcommand{\cc} [1]     {\overline{{#1}}}
\newcommand{\id}         {{\mathsf{id}}}
\newcommand{\st}[1]       {\scriptscriptstyle{{#1}}}

\newcommand{\cl}         {{\mathrm{cl}}}
\newcommand{\ring}[1]    {{\mathsf{{#1}}}}

\newcommand{\starWeyl} {\mathbin{\star_{\st{\mathrm{Weyl}}}}}
\newcommand{\starWick} {\mathbin{\star_{\st{\mathrm{Wick}}}}}

\newtheorem{lemma} {Lemma} [section]
\newtheorem{proposition} [lemma] {Proposition}

\newtheorem{theorem} [lemma] {Theorem}
\newtheorem{corollary} [lemma] {Corollary}

\theorembodyfont{\rm} 

\newtheorem{remark}[lemma]{Remark}

\newenvironment{proof}{{\sc Proof:}}{{\hspace*{\fill} $\square$\\}}

\numberwithin{equation}{section}

\begin{document}

\maketitle

\begin{abstract}
    Let $M$ be a Poisson manifold equipped with a Hermitian star
    product.  We show that any positive linear functional on
    $C^\infty(M)$ can be deformed into a positive linear functional
    with respect to the star product.
\end{abstract}

%
%

\section{Introduction}
\label{sec:intro}

A natural question in deformation quantization is whether
``classical'' positive linear functionals can be deformed into
``quantum'' positive linear functionals. In
\cite{bursztyn.waldmann:2000a} we gave an affirmative answer to this
question in the case of deformation quantization of symplectic
manifolds; in this paper, we extend this result to arbitrary Poisson
manifolds.

More precisely, let $M$ be a smooth manifold, and let $C^\infty(M)$ be
the algebra of complex-valued smooth functions on $M$. A
\textit{positive linear functional} on $C^\infty(M)$ is a complex
linear functional $\omega_0:C^\infty(M)\to \mathbb{C}$ satisfying
$\omega_0(\overline{f}\cdot f)\geq 0$ for all $f \in C^\infty(M)$.
(Such functionals are always given by integration with respect to
compactly supported positive Borel measures on $M$, see e.g.
\cite[App.~B]{bursztyn.waldmann:2001a}).  In the framework of
deformation quantization \cite{bayen.et.al:1978a}, $M$ is quantized by
a star product $\star$ on $C^\infty(M)[[\lambda]]$, the algebra of
formal power series in a real parameter $\lambda$ with coefficients in
$C^\infty(M)$. We assume in addition that $\overline{f\star
  g}=\overline{g}\star \overline{f}$ so that
$(C^\infty(M)[[\lambda]],\star)$ is an associative algebra over
$\mathbb{C}[[\lambda]]$ for which the pointwise complex conjugation of
functions is an involution. In order to define ``quantum'' positive
linear functionals, we use the natural notion of ``asymptotic
positivity'' in the ring $\mathbb{R}[[\lambda]]$: if
$a=\sum_{r=0}^\infty\lambda^r a_r \in \mathbb{R}[[\lambda]]$, then
$a>0$ if and only if $a_{r_0}>0$, where $a_{r_0}$ is the first
non-zero coefficient of $a$. Then, as before, a
$\mathbb{C}[[\lambda]]$-linear functional $\omega_0:
C^\infty(M)[[\lambda]]\to \mathbb{C}[[\lambda]]$ is called
\textit{positive} if $\omega_0(\overline{f}\star f)\geq 0$ for all
$f$.

If $\omega_0$ is a positive linear functional on $C^\infty(M)$, then
its $\lambda$-linear extension to $(C^\infty(M)[[\lambda]],\star)$
need not be positive. A concrete example is given when
$M=\mathbb{R}^{2n}$, $\omega_0$ is the delta functional at $0$, and
$\star$ is the Weyl-Moyal star product, see e.g.
\cite{bordemann.waldmann:1998a,bursztyn.waldmann:2000a}.  The natural
question is then whether one can find ``quantum corrections''
$\omega_i:C^\infty(M)\to \mathbb{C}$ so that $\omega=\omega_0 +
\sum_{i=1}^\infty \lambda^i \omega_i$ is a positive linear functional
on $(C^\infty(M)[[\lambda]],\star)$. A complete answer to this
question is provided by Theorem \ref{theorem:ComPosStars}, which
asserts that this is always possible.

The main ingredient in the proof of Theorem \ref{theorem:ComPosStars}
is the observation that any star product on $\mathbb{R}^n$ can be
realized as a subalgebra of the algebra of functions on the ``formal
cotangent bundle'' of $\mathbb{R}^n$ equipped with the Weyl-Moyal star
product, a result that relies on
\cite{bordemann.neumaier.nowak.waldmann:1997a:misc,nowak:1997a}.
Using this fact, the proof follows the same steps as the one for
symplectic star products in \cite[Prop.~5.1]{bursztyn.waldmann:2000a}.

The importance of positive linear functionals in deformation
quantization is illustrated by their central role in the
representation theory of star products initiated in
\cite{bordemann.waldmann:1998a}, see \cite{waldmann:2004a:pre} for a
recent review. In particular, Theorem \ref{theorem:ComPosStars} has
direct applications to the theory of strong Morita equivalence of star
products, see \cite{bursztyn.waldmann:2004}.

The paper is organized as follows: Section \ref{sec:def} contains the
basic definitions and the statement of the main theorem (Theorem
\ref{theorem:ComPosStars}); Section \ref{sec:Realization} contains the
main construction underlying its proof; Section \ref{sec:proof}
completes the proof of theorem.

\noindent
\textbf{Acknowledgements:} We thank Martin Bordemann for many helpful
discussions on the computation of Hochschild cohomologies using the
Koszul resolution. H.B. thanks DAAD for financial support and Freiburg
University for its hospitality while part of this work was being done.

%
%

\section{Basic definitions and the main theorem}
\label{sec:def}

Let us recall the general algebraic setting in which positive linear
functionals and positive deformations can be defined, see e.g.
\cite{bursztyn.waldmann:2000a}.

Let $\ring{C}$ be a ring of the form $\ring{R}(\im)$, where $\ring{R}$
is an ordered ring and $\im^2=-1$. Let $\mathcal{A}$ be an algebra
over $\ring{C}$ equipped with an involution $^*$. Using the order
structure on $\ring{R}$, we define a \textdef{positive linear
  functional} on $\mathcal{A}$ to be a $\ring{C}$-linear functional
$\omega:\mathcal{A}\to \ring C$ satisfying $ \omega(a^*\cdot a)\geq 0,
$ for all $a\in \mathcal{A}$.

If $\boldsymbol{\mathcal{A}}=(\mathcal{A}[[\lambda]],\star)$ is a
formal associative deformation of $\mathcal{A}$ in the sense of
Gerstenhaber \cite{gerstenhaber:1964a}, then we call it
\textdef{Hermitian} if
$$
(a_1\star a_2)^* =a_2^* \star a_1^*,\;\; \mbox{ for all }  a_1,a_2 \in \mathcal{A}.
$$
In this case, the $\lambda$-linear extension of the involution $^*$
from $\mathcal{A}$ to $\mathcal{A}[[\lambda]]$ makes
$\boldsymbol{\mathcal{A}}$ into a $^*$-algebra over
$\ring{C}[[\lambda]]$. Since $\ring{R}[[\lambda]]$ has an order
structure induced from that of $\ring{R}$ (analogous to the one in
$\mathbb{R}[[\lambda]]$ discussed in the introduction) and
$\ring{C}[[\lambda]]=\ring{R}[[\lambda]](\im)$, the definition of
positive linear functionals makes sense for $\boldsymbol{\mathcal{A}}$
as well.

Note that if 
$$
\omega = \sum_{r=0}^\infty \lambda^r \omega_r:
\mathcal{A}[[\lambda]] \longrightarrow \ring{C}[[\lambda]]
$$
is a positive $\ring{C}[[\lambda]]$-linear functional with respect
to $\star$, then its classical limit $\omega_0: \mathcal{A}
\longrightarrow \ring{C}$ is a positive $\ring{C}$-linear functional
on $\mathcal{A}.$ Conversely, we say that a Hermitian deformation
$\boldsymbol{\mathcal{A}} = (\mathcal{A}[[\lambda]], \star)$ is
\textdef{positive} \cite[Def.~4.1]{bursztyn.waldmann:2000a} if for
every positive linear functional $\omega_0$ of $\mathcal{A}$ one can
find $\ring{C}$-linear functionals
$$
\omega_r:\mathcal{A}\to \ring{C}, \qquad  r=1,2,\ldots,
$$
so that $\omega_0+\sum_{r=1}^\infty\lambda^r\omega_r$ is a positive
linear functional of $\boldsymbol{\mathcal{A}}$.  We say that
$\boldsymbol{\mathcal{A}}$ is a \textdef{completely positive}
deformation if, for each $n \in \mathbb{N}$, the $^*$-algebra
$M_n(\boldsymbol{\mathcal{A}})$ is a positive deformation of
$M_n(\mathcal{A})$.

A simple example illustrates that not all Hermitian deformations are
positive.  Let $\mathcal{A}$ be a $^*$-algebra over $\ring{C}$, and
let $\mu: \mathcal{A} \otimes \mathcal{A} \longrightarrow \mathcal{A}$
denote the multiplication map.  If we view $\mathcal{A}$ as an algebra
equipped with the \emph{zero} multiplication, then it is a
$^*$-algebra for which all linear functionals are positive. Then
$\lambda \mu$ provides a Hermitian deformation, which is clearly not
positive in general.

In this paper, we are concerned with algebraic deformations arising in
the geometric context of deformation quantization: If
$(M,\{\cdot,\cdot\})$ is a Poisson manifold, then a \textdef{star
  product} \cite{bayen.et.al:1978a} on $M$ is a formal associative
deformation $\star$ of the complex algebra $C^\infty(M)$,
$$
f\star g = f\cdot g + \sum_{r=1}^\infty \lambda^r C_r(f,g),
$$
for which each $C_r$ is a bidifferential operator on $M$ and
$$
C_1(f,g)-C_1(g,f)=\im\{f,g\}.
$$
We call the star product \textdef{Hermitian} if $\overline{f\star
  g}=\overline{g}\star \overline{f}$.

The following is the main result of this paper.
\begin{theorem}
    \label{theorem:ComPosStars}
    Any Hermitian star product on a Poisson manifold is a completely
    positive deformation.
\end{theorem}
The proof of Theorem \ref{theorem:ComPosStars} will be presented in
Section \ref{sec:proof}.

A key ingredient in the proof is the fact that any star product on
$\mathbb{R}^n$ can be realized as a subalgebra of the familiar
Weyl-Moyal star product on the ``formal cotangent bundle'' of
$\mathbb{R}^n$. More precisely, let us consider $\mathcal{W} =
\mathcal{W}_0[[\lambda]]$, where $\mathcal{W}_0 =
C^\infty(\mathbb{R}^n)[[p_1, \ldots, p_n]]$, equipped with the formal
Weyl-Moyal star product $\starWeyl$ defined by
\begin{equation}
    \label{eq:WeylMoyal}
    a \starWeyl b = \mu \circ
    \eu^{\frac{\im\lambda}{2}\sum_{k=1}^n\left(
          \partial_{q^k} \otimes \partial_{p_k} - \partial_{p_k}
          \otimes \partial_{q^k}\right)}
    a \otimes b, \qquad a, b \in
\mathcal{W}.
\end{equation}
Here $\mu(a\otimes b) = ab$ is the undeformed commutative product of
$\mathcal{W}$ and $q^1, \ldots, q^n$ are the canonical coordinates on
$\mathbb{R}^n$.  With respect to $\starWeyl$, $\mathcal{W}$ is an
associative $^*$-algebra over $\mathbb{C}[[\lambda]]$ with involution
given by complex conjugation. (We treat the formal parameters as
real.)

Let $\pi^*: C^\infty(\mathbb{R}^n) \longrightarrow \mathcal{W}_0$ be
the natural inclusion, which is clearly an algebra homomorphism 
(thought of as dual to the projection $\pi$ of the ``formal
cotangent bundle'' of $\mathbb{R}^n$ onto $\mathbb{R}^n$). Since
$\starWeyl$ is \textit{homogeneous} in the sense that the degree map
\begin{equation}
    \label{eq:DegreeMap}
    \mathrm{deg} =
    \sum\nolimits_i p_i \frac{\partial}{\partial p_i}
    + \lambda \frac{\partial}{\partial \lambda}
\end{equation}
is a derivation of $\starWeyl$
\cite{bordemann.neumaier.waldmann:1998a}, it follows in particular
that the $\lambda$-linear extension
$\pi^*:C^\infty(\mathbb{R}^n)[[\lambda]] \longrightarrow \mathcal{W}$
satisfies
$$
\pi^*(fg) = \pi^*(f) \starWeyl \pi^*(g).
$$
If $\star$ is a star product on $\mathbb{R}^n$ quantizing an
arbitrary Poisson structure $\{\cdot,\cdot\}$, then the claim is that
$\pi^*$ can be ``deformed'' into a $\mathbb{C}[[\lambda]]$-linear
algebra homomorphism
\begin{equation}
    \label{eq:TauDeformsPi}
    \tau: \left(C^\infty(\mathbb{R}^n)[[\lambda]], \star \right)
    \longrightarrow \left(\mathcal{W}, \starWeyl \right).
\end{equation}
\begin{remark}
    The \textit{classical limit} $\cl(\tau)$ of $\tau$ (defined by
    setting $\lambda$ to zero) gives an injective homomorphism of
    Poisson algebras
    \begin{equation}
        \label{eq:clTau}
        \cl(\tau): \left(C^\infty(\mathbb{R}^n), \{\cdot,\cdot\}\right)
        \longrightarrow
        \left(\mathcal{W}_0, \{\cdot,\cdot\}_{\st{can}} \right),
    \end{equation}
    where $\{\cdot,\cdot\}_{\st{can}}$ is the canonical Poisson
    bracket on $\mathcal{W}_0$. Since $\mathcal{W}_0$ can be thought
    of as the algebra of functions on the ``formal cotangent bundle''
    or $\mathbb{R}^n$, we think of \eqref{eq:clTau} as a ``formal
    symplectic realization'' of $(\mathbb{R}^n,\{\cdot,\cdot\})$, and
    of $\tau$ as its quantization.
\end{remark}

%
%

\section{Constructing the homomorphism $\tau$}
\label{sec:Realization}

The fact that there are no cohomological obstructions for the
recursive construction of the map $\tau$ as in \eqref{eq:TauDeformsPi}
relies on the joint work of Martin Bordemann, Nikolai Neumaier, Claus
Nowak and Stefan Waldmann
\cite{bordemann.neumaier.nowak.waldmann:1997a:misc}, see also the
thesis \cite{nowak:1997a}. We will outline the main ideas here.

Regarding $\mathcal{W}_0=C^\infty(\mathbb{R}^n)[[p_1,\ldots,p_n]]$ as
a $C^\infty(\mathbb{R}^n)$-bimodule via left and right multiplication
with respect to the usual pointwise product, we consider the
\textit{differential (resp. continuous)} Hochschild cohomology of
$C^\infty(\mathbb{R}^n)$ with values in $\mathcal{W}_0$, denoted by
$\mathrm{HH}^k (C^\infty(\mathbb{R}^n), \mathcal{W}_0)$, $k\geq 0$.
This is the cohomology of the complex
$$
\left(\oplus_{k=0}^\infty
\mathrm{HC}^k(C^\infty(\mathbb{R}^n), \mathcal{W}_0), \delta_0 \right),
$$
where $\mathrm{HC}^k(C^\infty(\mathbb{R}^n), \mathcal{W}_0)$ is the
space of $k$-multilinear maps from
$C^\infty(\mathbb{R}^n)\times\ldots\times C^\infty(\mathbb{R}^n)$ ($k$
times) into $\mathcal{W}_0$ which are differential operators on each
argument (resp. continuous with respect to the Fr\'echet structure of
$C^\infty(M)$), and $\delta_0$ is the Hochschild coboundary operator,
see e.g. \cite{gerstenhaber:1964a}.

Similarly, we regard $\mathcal{W}=\mathcal{W}_0[[\lambda]]$ as a
bimodule over $C^\infty(M)$ via left and right multiplication with
respect to $\starWeyl$. We denote the corresponding (differential,
resp.  continuous) Hochschild cochains by
$\mathrm{HC}^k(C^\infty(\mathbb{R}^n), \mathcal{W})$, the Hochschild
coboundary operator by $\delta$, and the Hochschild cohomology of
$C^\infty(M)$ with values in $\mathcal{W}$ by
$\mathrm{HH}^k(C^\infty(\mathbb{R}^n), \mathcal{W})$.  The continuous
cohomology (which turns out to coincide with the differential one) can
be explicitly described if one uses the (topological) Koszul
resolution of $C^\infty(\mathbb{R}^n)$, see e.g. \cite{connes:1994a},
and observes that the associated complex can be identified with
$$
\left(
    \oplus_{k=0}^\infty \mathcal{W}\otimes \wedge^k
    (\mathbb{R}^n)^*, \im\lambda d_p 
\right),
$$
where we view elements in $\mathcal{W}\otimes \bigwedge^k
(\mathbb{R}^n)^*$ as $k$-forms on $\mathbb{R}^{2n}$ of type
$\omega^{i_1\ldots i_k}(q,p)dp_{i_1}\wedge \ldots \wedge dp_{i_k}$,
with $\omega^{i_1\ldots i_k}(q,p) \in \mathcal{W}$, and $d_p$ is the
exterior derivative $d$ with respect to the $p_1, \ldots, p_n$
variables.  An application of the Poincar\'e's lemma shows the next
result.
\begin{proposition}
    \label{proposition:Hochschild}
    We have the following isomorphisms of
    $\mathbb{C}[[\lambda]]$-modules:
    \begin{equation}
        \label{eq:HHzero}
        \mathrm{HH}^0 (C^\infty(\mathbb{R}^n), \mathcal{W}) \cong
        C^\infty(\mathbb{R}^n)[[\lambda]],
    \end{equation}
    \begin{equation}
        \label{eq:Hochschildk}
        \mathrm{HH}^k (C^\infty(\mathbb{R}^n), \mathcal{W}) \cong
        (\mathcal{W}_0 \otimes \mbox{$\bigwedge$}^k (\mathbb{R}^n)^*)_{\mathrm{closed}}
        \quad
        \textrm{for}
        \quad
        k \ge 1.
    \end{equation}
\end{proposition}
Here $(\mathcal{W}_0 \otimes \mbox{$\bigwedge$}^k
(\mathbb{R}^n)^*)_{\mathrm{closed}}$ denotes the set of elements in
$\mathcal{W}\otimes \bigwedge^k (\mathbb{R}^n)^*$ which do not depend
on $\lambda$ and lie in the kernel of $d_p$.  The fact that continuous
and differential cohomologies coincide follows from techniques similar
to \cite{cahen.gutt.dewilde:1980a}, see also \cite{nowak:1997a}.

Observe that, due to the additional $\im\lambda$ multiplying $d_p$,
the cohomology is \emph{not} trivial (it would be zero if we used
formal Laurent series in $\lambda$ instead).  However,
\eqref{eq:Hochschildk} has the following consequence.
\begin{corollary}
    \label{cor:coboundary}
    If $\phi$ is a differential Hochschild $k$-cocycle with $k \ge 1$,
    then $\lambda\phi$ must be a coboundary.
\end{corollary}

Let $\cl: \mathcal{W} \longrightarrow \mathcal{W}_0$ denote the
classical limit map (setting $\lambda$ equal to zero), which is a
bimodule homomorphism with respect to the
$C^\infty(\mathbb{R}^n)$-bimodule structures. We keep the same
notation for the induced maps
\begin{equation}\label{eq:classlimit}
    \cl: \mathrm{HC}^{k }(C^\infty(\mathbb{R}^n),\mathcal{W}) \longrightarrow 
    \mathrm{HC}^{k }(C^\infty(\mathbb{R}^n),\mathcal{W}_0).
\end{equation}
Let $\mathrm{Alt}$ denote the antisymmetrization operator on
Hochschild cochains.
\begin{lemma}
    \label{lemma:NiceCoboundaries}
    Let $\phi \in \mathrm{HC}^{k \ge 1}(C^\infty(\mathbb{R}^n),
    \mathcal{W})$ be a cocycle with $\mathrm{Alt}(\cl(\phi)) = 0$.
    Then $\phi$ is a coboundary.
\end{lemma}
\begin{proof}
    Since $\cl: \mathcal{W} \longrightarrow \mathcal{W}_0$ is a
    bimodule homomorphism, it follows that the maps
    \eqref{eq:classlimit} satisfy
    \begin{equation}
        \label{eq:deltaclcldelta}
        \cl \circ \delta = \delta_0 \circ \cl.
    \end{equation}
    As a result, if $\phi \in \mathrm{HC}^k(C^\infty(\mathbb{R}^n),
    \mathcal{W})$ is a cocycle, i.e., $\delta \phi = 0$, then
    $\cl(\phi) \in \mathrm{HC}^k(C^\infty(\mathbb{R}^n),
    \mathcal{W}_0)$ satisfies $\delta_0 \cl(\phi) = 0$.  Since
    $\mathcal{W}_0 = C^\infty(\mathbb{R}^n) \otimes \mathbb{C}[[p_1,
    \ldots, p_n]]$, treating $p_1,\ldots,p_n$ as formal parameters we
    note that the usual Hochschild-Kostant-Rosenberg theorem for
    $C^\infty(\mathbb{R}^n)$ implies that, for $k\ge 1$, any cocycle
    in $\mathrm{HC}^k(C^\infty(\mathbb{R}^n), \mathcal{W}_0)$ is
    cohomologous to its skew symmetric part. Hence there exists a
    cochain $\psi \in \mathrm{HC}^{k-1}(C^\infty(\mathbb{R}^n),
    \mathcal{W}_0)$ such that
    \begin{equation}
        \label{eq:HKRcl}
        \cl(\phi) = \mathrm{Alt} (\cl(\phi)) + \delta_0 \psi.
    \end{equation}

    It follows that if $\mathrm{Alt}(\cl(\phi)) = 0$, then $\cl(\phi)
    = \delta_0 \psi$. Viewing $\psi$ as cochain with values in
    $\mathcal{W}$, we have that $\cl( \phi - \delta \psi) = 0$. Thus
    the cocycle $\phi - \delta \psi$ has the form $\lambda \eta$ for
    some other cocycle $\eta$. It follows from Corollary
    \ref{cor:coboundary} that $\phi$ is a coboundary.
\end{proof}
\begin{theorem}
    \label{theorem:SymplRealization}
    Let $\star$ be a star product on $\mathbb{R}^n$. For each $k\geq
    1$, there exists a $\tau_k \in
    \mathrm{HC}^1(C^\infty(\mathbb{R}^n), \mathcal{W})$, homogeneous
    of degree $k$ with respect to $\mathrm{deg}$, so that
    $$
    \tau= \pi^* + \sum_{k=1}^\infty\tau_k:
    \left(C^\infty(\mathbb{R}^n)[[\lambda]], \star \right)
    \longrightarrow \left(\mathcal{W}, \starWeyl \right)
    $$
    is an injective $\mathbb{C}[[\lambda]]$-linear algebra
    homomorphism.

    Furthermore, if $\star$ is a Hermitian star product, then one can
    chose $\tau$ to be a $^*$-homomorphism, i.e., $\tau(\cc{f}) =
    \cc{\tau(f)}$.
\end{theorem}
\begin{proof}
    For a sequence $\tau_i \in \mathrm{HC}^1(C^\infty(\mathbb{R}^n),
    \mathcal{W})$, $i=0,\ldots,k$, with each $\tau_i$ homogeneous of
    degree $i$ with respect to $\mathrm{deg}$, we define
    $\tau^{(k)}=\sum_{i=0}^k\tau_i$ and consider the error
    $\epsilon^{(k)}\in \mathrm{HC}^2(C^\infty(\mathbb{R}^n),
    \mathcal{W})$,
    \begin{equation}\label{eq:error}
        \epsilon^{(k)}(f,g)= \tau^{(k)}(f \star g) - \tau^{(k)}(f) \starWeyl \tau^{(k)}(g)
        \quad \textrm{for} \quad f,g \in C^\infty(\mathbb{R}^n).
    \end{equation}
    We write $\epsilon^{(k)}=\sum_{i=0}^\infty\epsilon^{(k)}_i$, with
    $\epsilon^{(k)}_i$ homogeneous of degree $i$ with respect to
    $\mathrm{deg}$.
    
    We now construct the desired $\tau$ recursively. If
    $\tau_0=\pi^*$, then $\epsilon^{(0)}_0=0$ since $\pi^*$ is a
    homomorphism for the undeformed products. Suppose that we have
    found $\tau_0, \ldots, \tau_{k-1}$ such that $\epsilon^{(k-1)}_0 =
    \cdots = \epsilon^{(k-1)}_{k-1} = 0$. Our goal is to find
    $\tau_k$, homogeneous of degree $k$, such that the error
    $\epsilon^{(k)}$ vanishes up to degree $k$.  Since
    $\epsilon^{(k)}_i= \epsilon^{(k-1)}_i$ for $i=0,\ldots,k-1$, we
    just have to impose the condition $\epsilon^{(k)}_{k} = 0$.  A
    direct computation shows that
    \begin{equation}
        \label{eq:Sk}
        \epsilon^{(k)}_k(f,g) = \sum_{i=0}^k \lambda^{k-i}\tau_i(C_{k-i}(f,g))
        - \sum_{i=1}^{k} \tau_i(f) \starWeyl \tau_{k-i}(g) =
        (\delta \tau_k) (f,g) + R_k(f,g),
    \end{equation}
    where $R_k$ depends only on $\tau_i$, $i=0,\ldots,k-1$. Explicitly,
    \begin{equation}
        \label{eq:Rk}
        R_k(f,g)
        = \sum_{i=0}^{k-1} \lambda^{k-i} \tau_i(C_{k-i}(f,g))
        - \sum_{i=1}^{k-1} \tau_i(f) \starWeyl \tau_{k-i}(g).
    \end{equation}
    Here $C_r$ is the $r$-th cochain of the star product $\star$.
    Hence we are left with showing that there exists $\tau_k$ of
    degree $k$ satisfying the cohomological equation
    \begin{equation}\label{eq;cobound}
        (\delta \tau_k) (f,g) + R_k(f,g)=0,
    \end{equation}
    i.e., that $R_k$ is a coboundary.  Note that, by \eqref{eq:error},
    $$
    \delta\epsilon^{(k)}(f,g,h)=f\starWeyl\epsilon^{(k)}(g,h)-\epsilon^{(k)}(f\star g, h) +
    \epsilon^{(k)}(f,g\star h)- \epsilon^{(k)}(f,g)\starWeyl h.
    $$
    Since $\delta$ and $\mathrm{deg}$ commute and
    $\epsilon^{(k)}_i=0$ for $i=0,\ldots,k-1$, we have
    $$
    \delta\epsilon^{(k)}_k(f,g,h)=f\starWeyl\epsilon^{(k)}_k(g,h)-\epsilon^{(k)}_k(f\star g, h) +
    \epsilon^{(k)}_k(f,g\star h)- \epsilon^{(k)}_k(f,g)\starWeyl h.
    $$
    Using \eqref{eq:error} and the associativity of $\star$ and
    $\starWeyl$, it is simple to check that
    $$
    \epsilon^{(k)}(f\star g, h) - \epsilon^{(k)}(f,g\star h) =
    \tau^{(k)}(f)\starWeyl\epsilon^{(k)}(g,h)-
    \epsilon^{(k)}(f,g)\starWeyl\tau^{(k)}(h),
    $$
    which in degree $k$ directly implies that
    $\delta\epsilon^{(k)}_k=0$. Hence, by \eqref{eq:Sk}, $\delta R_k =
    0$. Now a simple computation shows that $\cl(R_k)$ is symmetric,
    so Lemma~\ref{lemma:NiceCoboundaries} implies that $R_k$ is indeed
    exact. So we can find $\tau_k$ solving \eqref{eq;cobound}, and
    $\tau_k$ is homogeneous of degree $k$ because so is $R_k$.

    For the last part of the theorem, we must check that, if $\star$
    is Hermitian, then one can choose $\tau_k$ satisfying
    $\tau_k(\overline{f})=\overline{\tau_k(f)}$. Recall
    \cite[Sec.~3]{bursztyn.waldmann:2000a} that the complex
    conjugation on $C^\infty(\mathbb{R}^n)$ and $\mathcal{W}$ induce
    an involution on Hochschild cochains $\phi \in
    \mathrm{HC}^r(C^\infty(\mathbb{R}^n), \mathcal{W})$ by
    $$
    \phi^*(f_0,\ldots,f_r):=
    \overline{\phi(\overline{f_r},\ldots,\overline{f_0})},
    $$
    and $(\delta \phi)^*=(-1)^{r+1}\delta\phi^*$. If $\star$ is
    Hermitian and $\tau_i^*=\tau_i$ for $i=1,\ldots,k-1$, then a
    direct computation shows that $R_k$ defined in \eqref{eq:Rk} is
    Hermitian, i.e., $R_k^*=R_k$. If we pick $\tau_k$ with
    $\delta\tau_k = -R_k$, then $(\delta \tau_k)^*=\delta \tau_k^* =
    -R_k^*=-R_k$. It follows that
    $$
    \delta(\frac{1}{2}(\tau_k + \tau_k^*))=-R_k,
    $$
    so we can replace $\tau_k$ by its Hermitian part and assume
    that $\tau_k^*=\tau_k$.
\end{proof}
\begin{remark}
    \label{remark:Realization}
    The results in this section hold for any manifold $M$. This is
    obtained if one replaces $\mathcal{W}$ by the functions on $T^*M$
    depending formally on the `momentum variables' and $\starWeyl$ by
    any homogeneous star product on $T^*M$, see e.g.
    \cite{bordemann.neumaier.waldmann:1998a}.
\end{remark}


\section{Proof of the main theorem}
\label{sec:proof}
Let us consider complex coordinates $z^k = q^k + \im p_k$ and
$\cc{z}^k = q^k - \im p_k$ on $\mathbb{R}^{2n}$, and the derivations
$\frac{\partial}{\partial z^k}$ and $\frac{\partial}{\partial
  \cc{z}^k}$ of $\mathcal{W}_0[[\lambda]]$. For $A, B \in
M_N(\mathcal{W}_0)[[\lambda]]$, we consider the \emph{Wick star
  product}
\begin{equation}
    \label{eq:Wick}
    A \starWick B = \sum_{r=0}^\infty \frac{(2\lambda)^r}{r!}
    \sum_{i_1, \ldots, i_r}
    \frac{\partial^r A}
    {\partial z^{i_1} \cdots \partial z^{i_r}}
    \frac{\partial^r B}
    {\partial \cc{z}^{i_1} \cdots \partial \cc{z}^{i_r}}.
\end{equation}

We recall a general observation from
\cite[Sec.~4]{bursztyn.waldmann:2000a}.
\begin{lemma}
    \label{lem:Wick}
    If $\omega_0$ is a positive linear functional of
    $M_N(\mathcal{W}_0)$, then its $\lambda$-linear extension to
    $M_N(\mathcal{W}_0)[[\lambda]]$ is automatically positive with
    respect to $\starWick$.
\end{lemma}

We can now prove Theorem \ref{theorem:ComPosStars}.

\begin{proof}
    We first deal with the local case. Consider $\mathbb{R}^n$ (or any
    contractible open subset of it) equipped with an arbitrary Poisson
    structure, and let $\star$ be a Hermitian star product. Let us
    consider the formal Weyl algebra $(\mathcal{W} =
    \mathcal{W}_0[[\lambda]],\starWeyl)$.  By
    Theorem~\ref{theorem:SymplRealization}, we have an injective map $
    \tau: C^\infty(\mathbb{R}^n)[[\lambda]] \longrightarrow
    \mathcal{W}$, where $\tau = \sum_{k=0}^\infty \tau_k$, $\tau_0 =
    \pi^*: C^\infty(\mathbb{R}^n) \longrightarrow \mathcal{W}_0$ is
    the canonical inclusion, each $\tau_k$ is homogeneous of degree
    $k$, and
    \begin{equation}
        \label{eq:tauNiceMap}
        \tau(f \star g) = \tau(f) \starWeyl \tau(g)
        \quad
        \textrm{and}
        \quad
        \tau(\cc{f}) = \cc{\tau(f)}.
    \end{equation} 
    It is clear that $\tau$ extends to a $^*$-homomorphism for the
    matrix algebras,
    $$
    \tau: M_N(C^\infty(\mathbb{R}^n)[[\lambda]])
    \longrightarrow M_N(\mathcal{W}).
    $$

    Using the complex coordinates $z^k = q^k + \im p_k$ and $\cc{z}^k
    = q^k - \im p_k$, we recall that the operator
    \begin{equation}
        \label{eq:EquivWeylWick}
        S = \eu^{\lambda \Delta}
        \quad
        \textrm{with}
        \quad
        \Delta = \sum\nolimits_k
        \frac{\partial^2}{\partial z^k \partial \cc{z}^k}
    \end{equation}
    is an invertible $\mathbb{C}[[\lambda]]$-linear endomorphism of
    $M_N(\mathcal{W}_0)[[\lambda]]$ which is a $^*$-equivalence
    between $\starWick$ and $\starWeyl$
    \cite{bordemann.waldmann:1997a}, i.e.,
    \begin{equation}
        \label{eq:WickStarProduct}
        S(A^*)=S(A)^* \;\; \mbox{ and }\;\;  S (A \starWick B) = SA \starWeyl SB,
        \;\;
        \textrm{ for }
        \;
        A, B \in M_N(\mathcal{W}_0)[[\lambda]].
    \end{equation}

    Let $\Omega_0: M_N(C^\infty(\mathbb{R}^n)) \longrightarrow
    \mathbb{C}$ be a positive linear functional.  The canonical
    inclusion $\iota: \mathbb{R}^n \longrightarrow T^*\mathbb{R}^n$
    leads, at the algebra level, to a map $\iota^*: \mathcal{W}_0
    \longrightarrow C^\infty(\mathbb{R}^n)$ (just setting the
    `momentum variables' $p_1, \ldots, p_n$ to zero) which is a
    $^*$-homomorphism for the undeformed products and satisfies
    $\iota^*\pi^*=\id$. It follows that $\Omega_0 \circ \iota^*:
    M_N(\mathcal{W}_0) \longrightarrow \mathbb{C}$ is also a positive
    linear functional. Hence, by Lemma \ref{lem:Wick},
    $$
    \Omega_0 \circ \iota^*: M_N(\mathcal{W}_0)[[\lambda]]
    \longrightarrow \mathbb{C}[[\lambda]]
    $$
    is a positive linear functional with respect to $\starWick$,
    and by \eqref{eq:EquivWeylWick} we see that $\Omega_0 \circ
    \iota^* \circ S^{-1}$ is positive with respect to $\starWeyl$.
    Finally, by \eqref{eq:tauNiceMap}, the functional
    \begin{equation}
        \label{eq:omegaDefomeganull}
        \Omega = \Omega_0 \circ \iota^* \circ S^{-1} \circ \tau:
        C^\infty(\mathbb{R}^n)[[\lambda]]
        \longrightarrow \mathbb{C}[[\lambda]]
    \end{equation}
    is positive with respect to $\star$. Since $\iota^*\pi^* = \id$,
    it is easy to see that $\Omega$ is actually a deformation of
    $\Omega_0$. Therefore $\star$ is a completely positive
    deformation. Note that these results also hold for any
    contractible open subeset of $\mathbb{R}^n$.

    For the global result, we proceed just as in the symplectic case
    \cite[Prop.~5.1]{bursztyn.waldmann:2000a}. For a given Poisson
    manifold $M$ equipped with a Hermitian star product $\star$, we
    consider a locally finite open cover of $M$ by contractible open
    sets $\{\mathcal{O}_\alpha\}$ subordinate to a quadratic partition
    of unity $\{\chi_\alpha\}$ with $\sum_\alpha
    \overline{\chi_\alpha}\chi_\alpha =1$. If $\Omega_0$ is a positive
    linear functional on $M_N(C^\infty(M))$, then, on each
    $\mathcal{O}_\alpha$, we have a deformation $\Omega_{\alpha}$ of
    the restriciton of $\Omega_0$ to $\mathcal{O}_\alpha$. Then
    $$
    \Omega(f)=\sum_\alpha  \Omega_\alpha(\overline{\chi_\alpha}\star f \star \chi_\alpha)
    $$
    defines a positive linear functional on
    $(C^\infty(M)[[\lambda]],\star)$ which is a global deformation of
    $\Omega_0$.
\end{proof}

\begin{remark}
    One can directly prove Theorem \ref{theorem:ComPosStars} globally
    for $M$ by using Remark \ref{remark:Realization} together with a
    Wick star product defined via an almost complex structure on the
    ``formal cotangent bundle'' $T^*M$
    \cite{bordemann.waldmann:1997a}.
\end{remark}

%
%

\begin{footnotesize}

\begin{thebibliography}{10}

\bibitem {bayen.et.al:1978a}
{\sc Bayen, F., Flato, M., Fr{{\o}}nsdal, C., Lichnerowicz, A., Sternheimer,
  D.: }\newblock {\em Deformation Theory and Quantization}.
\newblock Ann. Phys.  {\bf 111} (1978), 61--151.

\bibitem {bordemann.neumaier.nowak.waldmann:1997a:misc}
{\sc Bordemann, M., Neumaier, N., Nowak, C., Waldmann, S.: }\newblock {\em
  Deformation of {P}oisson brackets}.
\newblock Unpublished discussions on the quantization problem of general
  Poisson brackets, June 1997.

\bibitem {bordemann.neumaier.waldmann:1998a}
{\sc Bordemann, M., Neumaier, N., Waldmann, S.: }\newblock {\em Homogeneous
  Fedosov Star Products on Cotangent Bundles I: Weyl and Standard Ordering with
  Differential Operator Representation}.
\newblock Commun. Math. Phys.  {\bf 198} (1998), 363--396.

\bibitem {bordemann.waldmann:1997a}
{\sc Bordemann, M., Waldmann, S.: }\newblock {\em A Fedosov Star Product of
  Wick Type for K{\"{a}}hler Manifolds}.
\newblock Lett. Math. Phys.  {\bf 41} (1997), 243--253.

\bibitem {bordemann.waldmann:1998a}
{\sc Bordemann, M., Waldmann, S.: }\newblock {\em Formal GNS Construction and
  States in Deformation Quantization}.
\newblock Commun. Math. Phys.  {\bf 195} (1998), 549--583.

\bibitem {bursztyn.waldmann:2000a}
{\sc Bursztyn, H., Waldmann, S.: }\newblock {\em On Positive Deformations of
  {$^*$}-Algebras}.
\newblock In: {\sc Dito, G., Sternheimer, D. (eds.): }\newblock {\em
  Conf{\'e}rence Mosh{\'e} Flato 1999. Quantization, Deformations, and
  Symmetries}, {\em Mathematical Physics Studies} no. {\bf 22},   69--80.
  Kluwer Academic Publishers, Dordrecht, Boston, London, 2000.

\bibitem {bursztyn.waldmann:2001a}
    {\sc Bursztyn, H., Waldmann, S.: }\newblock {\em Algebraic Rieffel Induction,
      Formal Morita Equivalence and Applications to Deformation Quantization}.
    \newblock J. Geom. Phys.  {\bf 37} (2001), 307--364.

\bibitem {bursztyn.waldmann:2004}
    {\sc Bursztyn, H., Waldmann, S.: }\newblock {\em Completely
      positive inner products and strong Morita equivalence}.
    \newblock Preprint {\bf math.QA/0309402} (2003),
    \newblock To appear in Pacific J. Math.

\bibitem{cahen.gutt.dewilde:1980a}
{\sc Cahen, M., Gutt, S., DeWilde, M.: }\newblock {\em Local Cohomology of the
  Algebra of $C^\infty$ Functions on a Connected Manifold}.
\newblock Lett. Math. Phys.  {\bf 4} (1980), 157--167.

\bibitem {connes:1994a}
    {\sc Connes, A.: }\newblock {\em Noncommutative Geometry}.
    \newblock Academic Press, San Diego, New York, London, 1994.

\bibitem {gerstenhaber:1964a}
    {\sc Gerstenhaber, M.: }\newblock {\em On the Deformation of Rings and
      Algebras}.
    \newblock Ann. Math.  {\bf 79} (1964), 59--103.
    
\bibitem{nowak:1997a}
    {\sc Nowak, C.~J.: }\newblock {\em {\"{U}}ber Sternprodukte auf
      nichtregul{\"{a}}ren Poissonmannigfaltigkeiten}.
    \newblock PhD thesis, Fakult{\"{a}}t f{\"{u}}r Physik,
    Albert-Ludwigs-Universit{\"{a}}t, Freiburg, 1997.

\bibitem{waldmann:2004a:pre}
    \textsc{Waldmann, S.}:
    \emph{States and representations in deformation quantization.}
    Preprint (Freiburg FR-THEP 2004/15)
    \textbf{math.QA/0408217} (August 2004), 51 pages.

\end{thebibliography}

\end{footnotesize}

\end{document}